\documentclass[%
]{ifacconf}
\usepackage[utf8]{inputenc}
\usepackage[T1]{fontenc}
\usepackage[final]{microtype}

\usepackage{tudcolors}
\usepackage{graphicx}
\usepackage{tikz}
\usetikzlibrary{
	plotmarks,
	math,
	positioning,
	shapes,
	arrows,
	arrows.meta,
	backgrounds,
	circuits.logic.IEC,circuits.ee.IEC,
	patterns,
	decorations,
	decorations.pathmorphing,
	shapes.geometric,
	calc,
	fit,
	matrix
}
\usepackage{gl}
\graphicspath{%
	{../include/}%
}

\usepackage{amsmath}
\usepackage{amssymb}
\usepackage{siunitx}
\usepackage{mathtools}

\makeatletter
\DeclareRobustCommand\vdots{%
  \mathpalette\@vdots{}%
}
\newcommand*{\@vdots}[2]{%
  \sbox0{$#1\cdotp\cdotp\cdotp\m@th$}%
  \sbox2{$#1.\m@th$}%
  \vbox{%
    \dimen@=\wd0 %
    \advance\dimen@ -3\ht2 %
    \kern.5\dimen@
    \dimen@=\wd2 %
    \advance\dimen@ -\ht2 %
    \dimen2=\wd0 %
    \advance\dimen2 -\dimen@
    \vbox to \dimen2{%
      \offinterlineskip
      \copy2 \vfill\copy2 \vfill\copy2 %
    }%
  }%
}
\DeclareRobustCommand\ddots{%
  \mathinner{%
    \mathpalette\@ddots{}%
    \mkern\thinmuskip
  }%
}
\newcommand*{\@ddots}[2]{%
  \sbox0{$#1\cdotp\cdotp\cdotp\m@th$}%
  \sbox2{$#1.\m@th$}%
  \vbox{%
    \dimen@=\wd0 %
    \advance\dimen@ -3\ht2 %
    \kern.5\dimen@
    \dimen@=\wd2 %
    \advance\dimen@ -\ht2 %
    \dimen2=\wd0 %
    \advance\dimen2 -\dimen@
    \vbox to \dimen2{%
      \offinterlineskip
      \hbox{$#1\mathpunct{.}\m@th$}%
      \vfill
      \hbox{$#1\mathpunct{\kern\wd2}\mathpunct{.}\m@th$}%
      \vfill
      \hbox{$#1\mathpunct{\kern\wd2}\mathpunct{\kern\wd2}\mathpunct{.}\m@th$}%
    }%
  }%
}
\makeatother

\DeclarePairedDelimiterX{\norm}[1]{\lVert}{\rVert}{#1}

\newcommand{\gen}{%
	\ensuremath{\circ}
}
\newcommand{\agen}{%
	\ensuremath{\bullet}
}

\newcommand{\dop}[1]{
	\,\mathrm{d}#1
}

\newcommand{\zbs}{%
	\ensuremath{\mathnormal{\Gamma}_{\mathrm{s}}}
}
\newcommand{\zbl}{%
	\ensuremath{\mathnormal{\Gamma}_{\mathrm{l}}}
}
\newcommand{\zbg}{%
	\ensuremath{\mathnormal{\Gamma}_{\gen}}
}

\newcommand{\zbt}{%
	\ensuremath{\tilde\Gamma}
}
\newcommand{\zbts}{%
	\ensuremath{\tilde{\mathnormal{\Gamma}}_{\mathrm{s}}}
}
\newcommand{\zbtl}{%
	\ensuremath{\tilde{\mathnormal{\Gamma}}_{\mathrm{l}}}
}
\newcommand{\zbtg}{%
	\ensuremath{\tilde{\mathnormal{\Gamma}}_{\gen}}
}

\newcommand{\errcomplete}[3]{%
	\ensuremath{\epsilon_{\mathrm{#1}}^{#2}%
		\ifthenelse{\equal{#3}{}}{}{(#3)}%
	}
}

\newcommand{\zisym}{%
	\ensuremath{\gamma}
}
\newcommand{\zicomplete}[4][]{%
	\ensuremath{%
			\zisym%
			_{\mathrm{#2}}^{#3}%
		\ifthenelse{\equal{#4}{}}{}{(#4)}%
	}
}
\newcommand{\vicomplete}[3]{%
	\ensuremath{\dot{\zisym}_{\mathrm{#1}}^{#2}%
		\ifthenelse{\equal{#3}{}}{}{(#3)}%
	}
}
\newcommand{\aicomplete}[3]{%
	\ensuremath{\ddot\zisym_{\mathrm{#1}}^{#2}%
		\ifthenelse{\equal{#3}{}}{}{(#3)}%
	}
}

\newcommand{\zi}[1][t]{%
	\ensuremath{\zicomplete{}{}{#1}}
}

\newcommand{\vi}[1][t]{%
	\ensuremath{\vicomplete{}{}{#1}}
}

\newcommand{\zir}{%
	\ensuremath{\zicomplete{r}{}{t}}
}
\newcommand{\vir}{%
	\ensuremath{\vicomplete{r}{}{t}}
}

\newcommand{\tm}{%
	\ensuremath{T_{\mathrm{m}}}%
}
\newcommand{\tnz}{%
	\ensuremath{T(z, t)}%
}
\newcommand{\torig}{%
	\ensuremath{T}%
}
\newcommand{\torigg}{%
	\ensuremath{T_{\gen}}%
}
\newcommand{\torigz}{%
	\ensuremath{\torig(z, t)}%
}

\newcommand{\ts}{%
	\ensuremath{T_{\mathrm{s}}}%
}

\newcommand{\tl}{%
	\ensuremath{T_{\mathrm{l}}}%
}

\newcommand{\tref}{%
	\ensuremath{T_{\mathrm{r}}}%
}

\newcommand{\trefg}{%
	\ensuremath{T_{\gen, \mathrm{r}}}%
}

\newcommand{\zt}{%
	\ensuremath{\tilde z}
}

\newcommand{\ttrans}{%
	\ensuremath{\tilde T}%
}
\newcommand{\ttransg}{%
	\ensuremath{\tilde T_{\gen}}%
}
\newcommand{\ttransag}{%
	\ensuremath{\tilde T_{\agen}}%
}
\newcommand{\ttranss}{%
	\ensuremath{\ttrans_\text{s}}%
}
\newcommand{\ttransl}{%
	\ensuremath{\ttrans_\text{l}}%
}
\newcommand{\ttransz}{%
	\ensuremath{\ttrans(\tilde z, t)}%
}
\newcommand{\ttransgz}{%
	\ensuremath{\ttransg(\tilde z, t)}%
}
\newcommand{\ttransref}{%
	\ensuremath{\tilde \tref}%
}
\newcommand{\ttransrefg}{%
	\ensuremath{\tilde T_{\gen, \mathrm{r}}}%
}

\newcommand{\ttransrefz}{%
	\ensuremath{\ttransref(\zt, t)}%
}

\newcommand{\ttransrefgz}{%
	\ensuremath{\ttransrefg(\zt , t)}%
}

\newcommand{\dzi}{%
	\ensuremath{\Delta\zi}%
}

\newcommand{\dvi}{%
	\ensuremath{\Delta\vi}%
}

\newcommand{\terr}{%
	\ensuremath{e}%
}

\newcommand{\terrz}{%
	\ensuremath{\terr(z, t)}%
}

\newcommand{\ttnerr}{e}
\newcommand{\ttnerrg}{%
	\ensuremath{\ttnerr_{\gen}}%
}
\newcommand{\ttnerrgz}{%
	\ensuremath{\ttnerrg(\zt, t)}%
}

\newcommand{\tterrg}{%
	\ensuremath{\tilde e_{\gen}}%
}
\newcommand{\tterrag}{%
	\ensuremath{\tilde e_{\agen}}%
}

\newcommand{\tterrgz}{%
	\ensuremath{\tterrg\left(\tilde z, t\right)}%
}

\newcommand{\hctrans}{%
	\mathnormal{\Psi}
}
\newcommand{\hctransi}{%
	\hctrans^{-1}
}
\newcommand{\hctransg}{%
	\hctrans_{\gen}
}
\newcommand{\hctransig}{%
	\hctransi_{\gen}
}

\newcommand{\hctransiag}{%
	\hctransi_{\agen}
}
\newcommand{\hctransgz}{%
	\hctransg\left(\zt, t\right)
}
\newcommand{\hctransigz}{%
	\hctransig\left(\zt, t\right)
}

\newcommand{\hctransiagz}{%
	\hctransiag\left(\zt, t\right)
}

\newcommand{\therrg}{%
	\ensuremath{\bar e_{\gen}}%
}
\newcommand{\therrag}{%
	\ensuremath{\bar e_{\agen}}%
}

\newcommand{\therrgz}{%
	\ensuremath{\therrg(\zt, t)}%
}

\newcommand{\ttarg}{%
	\ensuremath{w_{\gen}}%
}

\newcommand{\ttargz}{%
	\ensuremath{\ttarg(\zt, t)}%
}

\newcommand{\ttarmg}{%
	\ensuremath{\bar w_{\gen}}%
}
\newcommand{\ttarmag}{%
	\ensuremath{\bar w_{\agen}}%
}
\newcommand{\ttarmgz}{%
	\ensuremath{\ttarmg(\zt, t)}%
}

\newcommand{\hc}{%
	\ensuremath{\lambda}
}
\newcommand{\hcs}{%
	\ensuremath{\hc_{\mathrm{s}}}
}
\newcommand{\hcl}{%
	\ensuremath{\hc_{\mathrm{l}}}
}
\newcommand{\hcg}{%
	\ensuremath{\hc_{\gen}}
}
\newcommand{\hcag}{%
	\ensuremath{\hc_{\agen}}
}

\newcommand{\hd}{%
	\ensuremath{\alpha}
}
\newcommand{\hds}{%
	\ensuremath{\hd_{\mathrm{s}}}
}
\newcommand{\hdl}{%
	\ensuremath{\hd_{\mathrm{l}}}
}
\newcommand{\hdg}{%
	\ensuremath{\hd_{\gen}}
}

\newcommand{\dn}{%
	\ensuremath{\rho}
}
\newcommand{\ds}{%
	\ensuremath{\dn_{\mathrm{s}}}
}
\newcommand{\dl}{%
	\ensuremath{\dn_{\mathrm{l}}}
}
\newcommand{\dmelt}{%
	\ensuremath{\dn_{\mathrm{m}}}
}

\newcommand{\stc}{%
	\ensuremath{c_{\mathrm{p}}}
}
\newcommand{\stcs}{%
	\ensuremath{c_{\mathrm{p,s}}}
}
\newcommand{\stcl}{%
	\ensuremath{c_{\mathrm{p,l}}}
}

\newcommand{\bfd}{%
	\ensuremath{\delta_{\gen}}
}
\newcommand{\bfda}{%
	\ensuremath{\delta_{\agen}}
}
\newcommand{\bfds}{%
	\ensuremath{\delta_{\mathrm{s}}}
}
\newcommand{\bfdl}{%
	\ensuremath{\delta_{\mathrm{l}}}
}

\newcommand{\domsym}{\mathnormal{\Omega}}

\newcommand{\doms}{%
	\ensuremath{\domsym_{\mathrm{s}}}
}
\newcommand{\doml}{%
	\ensuremath{\domsym_{\mathrm{l}}}
}

\newcommand{\scg}{%
	\ensuremath{s_{\gen}}
}
\newcommand{\scag}{%
	\ensuremath{s_{\agen}}
}

\newcommand{\domt}{\tilde\Omega}

\newcommand{\inpg}{%
	\ensuremath{u_{\circ}(t)}
}
\newcommand{\inps}{%
	\ensuremath{u_{\mathrm{s}}(t)}
}
\newcommand{\inpl}{%
	\ensuremath{u_{\mathrm{l}}(t)}
}

\newcommand{\inptg}{%
	\ensuremath{\tilde u_{\circ}(t)}
}

\newcommand{\y}{%
	\ensuremath{\bs{y}(t)}
}

\newcommand{\bs}[1]{\boldsymbol{#1}}
\newcommand{\partiell}[3][]{\frac{\partial^{#1}#2}{\partial{#3}^{#1}}}

\newcommand{\cgt}[1]{%
	\ensuremath{c_{\gen,#1}(t)}
}

\newcommand{\kernelcomplete}[4]{%
    \ensuremath{#1(#2, #3, #4)}
}

\newcommand{\kt}{%
    \tilde k
}
\newcommand{\ktzzt}{%
    \kernelcomplete{\kt}{\zt}{\zeta}{t}
}
\newcommand{\ktg}{%
    \tilde k_{\gen}
}
\newcommand{\ktgzzt}{%
    \kernelcomplete{\ktg}{\zt}{\zeta}{t}
}

\newcommand{\kb}{%
	\bar k
}
\newcommand{\kbest}{%
	\kernelcomplete{\kb}{\eta}{\sigma}{t}
}
\newcommand{\kbeot}{%
	\kernelcomplete{\kb}{\eta}{0}{t}
}
\newcommand{\kbeet}{%
	\kernelcomplete{\kb}{\eta}{\eta}{t}
}

\newcommand{\Kb}{%
	\bar K
}
\newcommand{\Kernelcomplete}[5]{%
	\ensuremath{#1^{#5}(#2, #3, #4)}
}
\newcommand{\Kbest}[1]{%
	\Kernelcomplete{\Kb}{\eta}{\sigma}{t}{#1}
}

\newcommand{\discrs}{%
	\Delta_\sigma
}
\newcommand{\discre}{%
	\Delta_\eta
}
\newcommand{\discr}{%
	\Delta
}
\newcommand{\kij}[2]{%
	\bar k_{#1, #2}(t)
}

\newcommand{\gorder}{%
	\aleph
}
\newcommand{\gevclasscomplete}[2]{%
	\mathcal{G}_{#1}(#2)
}
\newcommand{\ggen}{%
	\gevclasscomplete{\gorder}{\domsym}
}
\newcommand{\fclasscomplete}[2]{%
	\mathcal{K}_{#1}(#2)
}

\def\highlight<#1>{%
	\temporal<#1>{\color{HKS41}}{\color{HKS44}}{\color{HKS41}}
}

\RequirePackage[babel]{csquotes} 
\usepackage{natbib}        

\usepackage{booktabs}

\usepackage{glossaries}
\newacronym{dps}{DPS}{distributed parameter system}
\newacronym{bvp}{BVP}{boundary value problem}
\newacronym{ivp}{IVP}{initial value problem}
\newacronym{fbp}{FBP}{free boundary problem}
\newacronym{sp}{SP}{Stefan problem}
\newacronym{opsp}{OPSP}{one-phase Stefan problem}
\newacronym{tpsp}{TPSP}{two-phase Stefan problem}
\newacronym{sc}{SC}{Stefan condition}
\newacronym{ode}{ODE}{ordinary differential equation}
\newacronym{odes}{ODES}{ordinary differential equation system}
\newacronym{pde}{PDE}{partial differential equation}
\newacronym{fem}{FEM}{finite element method}
\newacronym{cfd}{CFD}{computational fluid dynamic}
\newacronym{vgf}{VGF}{Vertical Gradient Freeze}
\newacronym{tmf}{TMF}{travelling magnetic field}
\newacronym{gaas}{GaAs}{Gallium-Arsenide}
\newacronym{inp}{InP}{Indium-Phosphide}
\newacronym{dof}{DOF}{degree of freedom}
\newacronym{iss}{ISS}{Input State Stability}
\newacronym{pi}{CSI}{Poincaré inequality}
\newacronym{csi}{CSI}{Cauchy-Schwarz inequality}
\newacronym{yi}{YI}{Young's inequality}
\newacronym{lti}{LTI}{Linear Time Invartiant}
\newacronym{ltv}{LTV}{Linear Time Vartiant}
\newacronym{lqe}{LQE}{Linear Quadratic Estimator}
\newacronym{frde}{FDRE}{Filtering Riccati Differential Equation}

\begin{document}
\begin{frontmatter}

\title{
Control of the Vertical Gradient Freeze crystal growth process via backstepping}

\thanks[footnoteinfo]{
This work has been funded by the Deutsche Forschungsgemeinschaft (DFG) [project number WI 4412/1-1
].}
\thanks[footnoteinfo]{
	©2020 the authors.
	This work has been accepted to IFAC for publication under a Creative Commons Licence CC-BY-NC-ND
}

\author[First]{Stefan Ecklebe} 
\author[Second]{Frank Woittennek} 
\author[First]{Jan Winkler}

\address[First]{Institute of Control Theory, Technische Universit\"at Dresden,
01069 Dresden, Germany (email: \{stefan.ecklebe, jan.winkler\}@tu-dresden.de)}
\address[Second]{Institute of Automation and Control Engineering, University of
Health Sciences, Medical Informatics and Technology, 6060 Hall in Tirol, Austria (email: frank.woittennek@umit.at)}

\begin{abstract}    
	This contribution presents a backstepping-based state feedback design for
	the tracking control of a \acrlong{tpsp} which is encountered in the 
	\acrlong{vgf} crystal growth process.
	A \acrlong{tpsp} consists of two coupled \acrlongpl{fbp} and
	is a vital part of many crystal growth processes due to the time-varying
	extent of crystal and melt during growth.
	In addition, a different approach for the numerical approximation of
	the backstepping transformations kernel is presented.
\end{abstract}


\end{frontmatter}
\section{Introduction}

The \gls{vgf} crystal growth process is used for the production of high efficiency bulk compound
semiconductor single crystals like \gls{gaas} or Indium-Phosphide \citep{JURISCH2005283}.
The process basically works as follows:
A seed crystal is placed at the bottom of a rotationally symmetric crucible which is later filled with
solid semiconductor chunks.
After all material (up to the seed) in the crucible is molten,
a vertical temperature gradient is moved through the plant such that a single crystal grows from the
bottom to the top in a desired manner.
This is done by manipulating the power of the heaters which surround the crucible.
Modelling of the system yields two coupled  \glspl{fbp} for crystal and melt that form
the so called \gls{tpsp} e.g.~\citep{Crank84} which is inherently nonlinear.

Due to the spatial extent of the system it is broadly discussed in the framework of distributed
parameter systems.
Making the assumption that the temperature distribution in one phase is constant (which is often
justified due to its dominant spatial extent) yields the so called \gls{opsp}.
Regarding this special case, results are lately available for the feedforward design 
by \cite{DPRM03}, as well as for feedback designs using enthalpy-
\citep{Petrus2012, Petrus2014}, geometry- \citep{Maidi2016} or backstepping- \citep{Koga2019a}
based approaches.
Regarding the full problem, \citep{RWW03md, RWW04esta} extend the flatness-based motion planning to the
two-phase case, while \citep{Kang1995} and \citep{Hinze2009} address the problem from the side
of optimal control.
Concerning feedback, a direct extension of the approaches for the one-phase variant is not feasible since 
the coupling between the two \glspl{fbp} has to be taken into account here.
In this context it is noteworthy that \citep{petrus2010} already states a Lyapunov-based control 
law for the \gls{tpsp} with actuation at one boundary and that \citep{Koga2019c} does the same
via an energy-shaping approach.
However, a gap to the tracking control of the complete \gls{tpsp} using both inputs remained.
In \citep{Ecklebe2019} the authors present different output feedback designs using energy- and
flatness-based approaches, 
rendering backstepping-based designs the last remaining problem.

\subsection{Objective and Structure} 
The main objective of this contribution is to present a tracking control for
the \gls{tpsp} by backstepping-based state feedback. 
Furthermore, a different approach for the numeric approximation of the resulting time-variant
backstepping kernels is presented.
Due to the limited space, this is done in a rather brief fashion and more detailed results will be
given in a forthcoming publication.

In Section \ref{sec:modelling} a simplified one-dimensional distributed parameter model of the 
process is introduced, before Section \ref{sec:feedforward} briefly recites the feedforward control
design and states some properties of the derived trajectories which are required for the tracking
control later on.
Based on these, Section \ref{sec:backstepping} derives a suitable error system as well as the
corresponding target dynamics for it and states the resulting tracking control law.
In Section \ref{sec:kernel}, the existence of solutions for the kernel equations as well as a new
computation scheme to solve them is discussed, before Section \ref{sec:results} briefly presents
simulation results.
Finally, a summary and an outlook to further work is given.


\section{Modelling}
\label{sec:modelling}

\begin{figure}
	\centering
	\resizebox{\linewidth}{!}{\tikzset{>=stealth}%
\newlength{\hatchspread}%
\newlength{\hatchthickness}%
\newlength{\hatchshift}%
\newcommand{\hatchcolor}{}%
\tikzset{hatchspread/.code={\setlength{\hatchspread}{#1}},
         hatchthickness/.code={\setlength{\hatchthickness}{#1}},
         hatchshift/.code={\setlength{\hatchshift}{#1}},
         hatchcolor/.code={\renewcommand{\hatchcolor}{#1}}}%
\tikzset{hatchspread=15pt,
         hatchthickness=5pt,
         hatchshift=0pt,
         hatchcolor=HKS44K100}%
\pgfdeclarepatternformonly[\hatchspread,\hatchthickness,\hatchshift,\hatchcolor]
   {custom north east lines}
   {\pgfqpoint{\dimexpr-2\hatchthickness}{\dimexpr-2\hatchthickness}}
   {\pgfqpoint{\dimexpr\hatchspread+2\hatchthickness}{\dimexpr\hatchspread+2\hatchthickness}}
   {\pgfqpoint{\dimexpr\hatchspread}{\dimexpr\hatchspread}}
   {
    \pgfsetlinewidth{\hatchthickness}
    \pgfpathmoveto{\pgfqpoint{\dimexpr\hatchshift-0.15pt}{-0.15pt}}
    \pgfpathlineto{\pgfqpoint{\dimexpr\hatchspread+0.15pt}{\dimexpr\hatchspread-\hatchshift+0.15pt}}
    \ifdim \hatchshift > 0pt
      \pgfpathmoveto{\pgfqpoint{-0.15pt}{\dimexpr\hatchspread-\hatchshift-0.15pt}}
      \pgfpathlineto{\pgfqpoint{\dimexpr\hatchshift+0.15pt}{\dimexpr\hatchspread+0.15pt}}
    \fi
    \pgfsetstrokecolor{\hatchcolor}
    \pgfusepath{stroke}
   }%
\begin{tikzpicture}
	\glloadidentity
	\glrotatex{-70}

	\pgfmathsetmacro{\radius}{2}
	\pgfmathsetmacro{\height}{1.61*1.5*\radius}
	\pgfmathsetmacro{\phase}{\height*(1-1/1.61)}
	\pgfmathsetmacro{\arrow}{3}


	\glcylinderback[draw, dashed, thick, fill=HKS41K100, fill opacity=0]{\radius}{0}{\height}
	\glcylinderback[draw, dashed, thick, fill=HKS41K100, fill opacity=0]{\radius}{0}{\phase}
	\fill[glcoords, fill=HKS41K100, fill opacity=.1](0,0)circle(\radius);


	\begin{glscope}
		\gltranslatez{\phase}
		\fill[glcoords, pattern=custom north east lines, pattern color=HKS44K100,
		](0,0)circle(\radius);
	\end{glscope}

	
	\glcylinderfront[draw, thick, fill=HKS41K100, fill opacity=0]{\radius}{0}{\phase}
	
	\glcylinderfront[draw, thick, fill=HKS41K100, fill opacity=0]{\radius}{0}{\height}

	\begin{glscope}
		\gltranslatez{\height}
		\fill[glcoords, draw, thick, fill=HKS41K100, fill opacity=0.1](0, 0)circle(\radius);
	\end{glscope}

	\draw[glcoords, ->, thick](0,0) -- (\radius, 0) node[below, pos=.5]{$r$};
	\draw[glcoords, ->, thick](0:\radius/1.61) arc (0:90:\radius/2) node[left]{$\varphi$};
	\begin{glscope}
		\glrotatex{90}
		\gltranslatex{1.2*\radius}
		\draw[glcoords, {Rays[n=2]}-, thick] 
			(0, 0) node[right, pos=0]{$\zbs$} -- (0, \phase) node[right]{$\zi$};
		\draw[glcoords, {Rays[n=2]}-, thick] 
			(0, \phase) -- (0, \height) node[right]{$\zbl$};
		\draw[glcoords, {Rays[n=2]}->, thick] 
			(0, \height) -- (0, 1.1*\height) node[right]{$z$} ;

		\gltranslatex{.5*\radius}
		\draw[glcoords, {Rays[n=2]}-, thick] 
			(0, 0) node[right, pos=0]{$\zbts$} -- (0, \phase) node[right]{$0$};
		\draw[glcoords, {Rays[n=2]}-, thick] 
			(0, \phase) -- (0, \height) node[right]{$\zbtl$};
		\draw[glcoords, {Rays[n=2]}->, thick] 
			(0, \height) -- (0, 1.1*\height) node[right]{$\zt$} ;
	\end{glscope}

	\begin{glscope}
		\glrotatex{90}
		\gltranslatex{-2.6*\radius}
		\gltranslatey{.1*\height}
		\draw[glcoords, thick] (0, 0) -- node[above, pos=0, anchor=south west]{\large Bottom Heater}
			++(\arrow, 0) -- ++(.5\arrow, -.5\arrow);
		\gltranslatey{.2*\height}
		\draw[glcoords, thick] (0, 0) -- node[above, pos=0, anchor=south west]{\large Crystal}
			++(\arrow, 0) -- ++(.5\arrow, -.5\arrow);
		\gltranslatey{.2*\height}
		\draw[glcoords, thick] (0, 0) -- node[above, pos=0, anchor=south west]{\large Interface}
			++(\arrow, 0) -- ++(.5\arrow, -.5\arrow);
		\gltranslatey{.2*\height}
		\draw[glcoords, thick] (0, 0) -- node[above, pos=0, anchor=south west]{\large Melt}
			++(\arrow, 0) -- ++(.5\arrow, -.5\arrow);
		\gltranslatey{.2*\height}
		\draw[glcoords, thick] (0, 0) -- node[above, pos=0, anchor=south west]{\large Top Heater}
			++(\arrow, 0) -- ++(.5\arrow, .5\arrow);
	\end{glscope}

\end{tikzpicture}%
}
	\caption{Schematics of the cylindrical coordinates $(r, \varphi, z)$
		and the moving coordinate frame $\zt = z - \zi$.}
	\label{fig:system}
\end{figure}
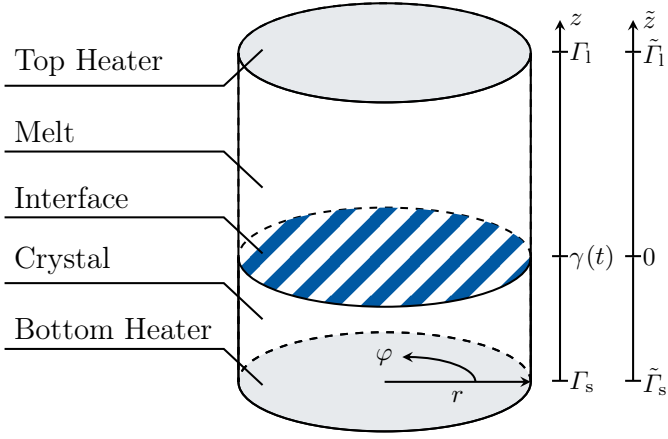

As the foundation for model based control, this section introduces a one dimensional distributed
parameter model of the \gls{vgf} process plant.

\subsection{Plant Model}
The quantity under consideration is given by the spatial and temporal distribution of the system
temperature $T$ in the crucible, 
denoted in cylindrical coordinates with radius $r$, angle $\varphi$ as well as height $z$, and
depending on the time $t$.
Within this contribution, we assume that the lateral heaters are used as active isolation,
avoiding any heat loss in radial direction and therefore yielding a temperature distribution which is
independent of $r$. 
Since the plant is also rotationally symmetric, this justifies averaging over the longitudinal 
cross-sectional area, reducing the spatial domain to a line whose boundaries are represented by the bottom
and top of the crucible at $z=\zbs$ and $z=\zbl$, respectively.
This yields two areas given by the crystal and the melt, separated by the moving phase
boundary $\zi$ (cf.~Figure \ref{fig:system}).
In contrast to the temperature distribution in the crystal, which can be modelled via diffusion,
the liquid melt also enables convective heat transport.
However, since the considered semi-conductors posses small Prandtl numbers (e.g.~\SI{0.068}{} for \gls{gaas}),
the dominating heat transport mechanism is diffusion. 
Therefore, convective effects in the melt are neglected.

Summarising, the temperature distribution in the system is given by the distributed variable
$T(z,t)$ and governed by a one dimensional nonlinear heat equation \citep{Rota84}
\begin{multline}
		\partiell{}{t}\Big(\rho(\tnz)\stc(\tnz) \tnz \Big)=\\
		\qquad \partiell{}{z}\left(\hc(\tnz) \partiell{}{z} \tnz \right)
		,z \in (\zbs, \zbl) \setminus \left\{\zi\right\}
	\label{eq:nonlin_heat_eq}
\end{multline}
with the density $\rho$, the specific heat capacity $\stc$, and $\hc$ the thermal conductivity
being temperature-dependent.
Note, that while the temperature at the interface $T(\zi,t)$ is fixed at the melting point temperature
$\tm$ due to the ongoing phase transition,
the heat flow in this description is not continuous at the phase boundary due to the release 
of latent heat within the solidification process and the abrupt change of the physical
parameters between crystal and melt.

\subsection{Decomposition}
Assuming piecewise constant parameters for the solid and the liquid phase it is possible to 
decompose the tempearature distribution via
\begin{equation}
	\tnz = \begin{cases}
		\ts(z,t), &z \in \doms = (\zbs, \zi)\\
		\tl(z,t), &z \in \doml = (\zi, \zbl)
	\end{cases}
	\label{eq:sep_temps}
\end{equation}
with the temperatures $\ts(z,t)$ and $\tl(z,t)$ in the solid and liquid part, respectively.
This yields the two \glspl{fbp} 
\begin{subequations}
	\label{eq:lin_heat_eqs}
	\begin{align}
		\label{eq:lin_heat_eq_s_pde}
		\partial_t \ts(z, t) &= \hds\partial_z^2 \ts(z,t)
		\\
		\label{eq:lin_heat_eq_s_bc1}
		\partial_z \ts (\zbs, t) &= \frac{\bfds}{\hcs}\inps\\
		\label{eq:lin_heat_eq_s_bc2}
		\ts(\zi, t) &= \tm \\
		\label{eq:lin_heat_eq_l_pde}
		\partial_t \tl(z, t) &= \hdl\partial_z^2 \tl(z,t)
		\\
		\label{eq:lin_heat_eq_l_bc1}
		\partial_z \tl (\zbl, t) &= \frac{\bfdl}{\hcl}\inpl\\
		\label{eq:lin_heat_eq_l_bc2}
		\tl(\zi, t) &= \tm
	\end{align}
\end{subequations} 
where the indexes ``s'' and ``l'' denote the solid and liquid phase, respectively.
Furthermore, the heat flows $\inps$ and $\inpl$ at the bottom and the top boundary are
considered as system inputs 
with the orientation factors $\bfds=-1$ and $\bfdl=1$. 
For reasons of clarity, the partial derivative of $\torigz$ wrt.~$z$ and $t$ are given by
$\partial_z \tnz$ and $\partial_t \tnz$, respectively.
Finally, $\hds = \hcs /(\ds\stcs)$ and $\hdl = \hcl/(\dl\stcl)$ denote the thermal 
diffusivities.

Next, examining the energy balance at the interface $\zi$ yields the \acrlong{sc} \citep{stefan_cond}
\begin{equation}
	\dmelt L\vi = \hcs\partial_z\ts(\zi, t) - \hcl\partial_z\tl(\zi, t)
	\label{eq:stefan_cond}
\end{equation}
which describes the evolution of the phase boundary.
Herein, $\dmelt$ denotes the density of the melt at melting temperature and $L$ the specific latent
heat.

Together, \eqref{eq:lin_heat_eqs} and \eqref{eq:stefan_cond} form the \gls{tpsp} whose state is 
given by
\begin{equation}
	\bs{x}(\cdot, t) = \begin{pmatrix}
		\torig(\cdot, t) \\ \zi
	\end{pmatrix} \in X = 
		L_2([\zbs, \zbl])
		\times (\zbs, \zbl) \:,
	\label{eq:state}
      \end{equation}
      where $L_2([\zbs, \zbl])$ denotes the space of real-valued square-integrable functions defined on $[\zbs, \zbl]$.
Note that the PDE-ODE-PDE system defined by 
\eqref{eq:lin_heat_eqs} and \eqref{eq:stefan_cond} is inherently nonlinear since the domains of
\eqref{eq:lin_heat_eq_s_pde} and \eqref{eq:lin_heat_eq_l_pde} depend on the state variable $\zi$.


Since the systems \eqref{eq:lin_heat_eq_s_pde}--\eqref{eq:lin_heat_eq_s_bc2}
and \eqref{eq:lin_heat_eq_l_pde}--\eqref{eq:lin_heat_eq_l_bc2} share the same structure,
the following sections will merely discuss generic variables, denoted by the $\gen$ symbol 
if the results are applicable to both phases.
If terms from two different phases are to appear in the same expression,
the complementary phase is marked by the $\agen$ symbol.

\subsection{Moving Coordinates}
To simplify the notation of the controller error system later on, the coordinate transformation
\begin{equation}
	\ttransz = \torigz \quad \text{with} \quad \zt \coloneqq z - \zi
	\label{eq:coord_trans}
\end{equation}
is introduced which maps the current interface position to the origin of a moving frame
as shown in Figure \ref{fig:system}.
This yields the generic system in the new coordinates
\begin{subequations}
	\begin{align}
		\label{eq:transformed_pde}
		\partial_t \ttransgz &= \hdg\partial_{\zt}^2\ttransgz + \vi \partial_{\zt} \ttransgz \\
		\label{eq:transformed_bc_0}
		\partial_{\zt} \ttransg(\zbtg, t) &= \frac{\bfd}{\hcg}\inpg \\
		\label{eq:transformed_bc_1}
		\ttransg(0, t) &= \tm \\
		\label{eq:transformed_sc}
		\vi &= \scg\partial_{\zt}\ttransg(0, t) +\scag \partial_{\zt}\ttransag(0, t)
	\end{align}%
	\label{eq:transformed_system}%
\end{subequations}
where $\zbtg = \zbg - \zi$, $\scg = -\bfd\hcg/(L \dmelt)$,
and $\scag = -\bfda\hcag/(L \dmelt)$.
Note, that in these coordinates the interface velocity $\vi$ directly enters the PDE 
\eqref{eq:transformed_pde} in form of a convection coefficient.


\section{Feedforward control}
\label{sec:feedforward}

This section outlines a feedforward control that originates from \cite{DPRM03} for
the \gls{opsp} and was extended by \cite{RWW03md} to the \gls{tpsp}.
Since the trajectories for $\torigz$ and $\zi$, which are computed by this feedforward scheme,
will be used as a reference in following sections, this recap is merely focussed on their
properties.

Since the solution $\ttransgz$ of \eqref{eq:transformed_system} can be expressed in terms of an
infinite power series
in $\zt$, given by
\begin{equation}
	\ttransgz= \sum\limits_{i=0}^\infty \cgt{i}\frac{\zt^i}{i!} \:,
	\label{eq:power_series}
\end{equation}
substitution into \eqref{eq:transformed_pde} and comparison of the coefficients of like powers
in $\zt$ yields the recursion formula
\begin{equation}
	\cgt{i+2} = \frac{1}{\hdg}\big(
		\partial_t \cgt{i}
		- \vi \cgt{i+1}
	\big)
	\ i=0,\dotsc,\infty \,.
	\label{eq:coeff_iteration}
\end{equation}
Examining \eqref{eq:power_series} indicates that the initial coefficients of the series are given by
\begin{alignat}{2}%
	\cgt{0} &= \ttransg(0, t) = \tm ,
	&\quad \cgt{1} &= \partial_{\zt} \ttransg(0, t) \:.
	\label{eq:c1}
\end{alignat}%
Next, using \eqref{eq:transformed_sc} as the defining equation for the melt-gradient
\begin{equation}
\partial_{\zt}\ttransl(0, t) = \frac{1}{\hcl}\left(\hcs\partial_{\zt}\ttranss(0, t) 
- \dmelt L\vi\right) \;,
\end{equation}
the solution for both phases can be expressed by the gradient in the crystal
$\partial_{\zt}\ttranss(0, t)$ and the growth rate $\vi$.
Thus, \eqref{eq:transformed_system} via the parametrisation \eqref{eq:power_series} is differentially
flat with a flat output
\begin{equation}
	\bs{y}(t) 
	= \begin{pmatrix} y_1(t) \\ y_2(t)\end{pmatrix}
	= \begin{pmatrix}
		\partial_{\zt}\ttranss(0, t) \\
		\zi
	\end{pmatrix}.
	\label{eq:flat_output}
\end{equation}
Furthermore, reference trajectories for the components $y_{\mathrm{r},i}(t)$ of $\bs{y}_{\mathrm{r}}(t)$ are chosen as
transitions between stationary states $y_{\mathrm{r}, i}^0$ and $y_{\mathrm{r}, i}^{\mathrm{e}}$ for start and end, respectively, via
\begin{equation}
	y_{\mathrm{r}, i}(t) = y_{\mathrm{r}, i}^0 +(y_{\mathrm{r}, i}^{\mathrm{e}}
		- y_{\mathrm{r}, i}^0)\Phi\left(t\right) \quad i=1,2
	\label{eq:flat_traj}
\end{equation}
with $\Phi(t)$ sufficiently smooth.

Analysing the specific convergence conditions for \eqref{eq:power_series},
\citep{RWW03md, RWW04esta} show that it is sufficient to demand $\Phi(t) \in \gevclasscomplete{\aleph\le2}{\mathbb{R}^+}$
with the Gevrey class $\ggen$ from the definition below.

\begin{defn}[Gevrey-class; \citealt{Gevrey1918}]
	A smooth function $t \mapsto f(t)$ defined on the open set 
	$\mathnormal{\Omega}\subset\mathbb{R}$ is an element of the 
	Gevrey class $\ggen$ of order $\gorder$ over $\mathnormal{\Omega}$
	if there exists a positive constant $D$ such that 
	\begin{equation*}
		\sup_{t \in \mathnormal{\Omega}} \lvert \partial_t^n f(t)\rvert \le D^{n+1}\left(n!\right)^\gorder
	\end{equation*}%
	holds for all $n $ in $\mathbb{N}_0$.
\end{defn}%

Thus, given the definition of the flat output \eqref{eq:flat_output}
the reference interface trajectory $\zir$ is from $\gevclasscomplete{\aleph\le2}{\mathbb{R}^+}$.
Moreover, as the reference temperature distribution $\ttransrefz$ is computed via \eqref{eq:coeff_iteration} and \eqref{eq:power_series},
by construction $\ttransrefz$ belongs to $\fclasscomplete{\aleph\le2}{\domt, \mathbb{R}^+}$
as defined below.
\begin{defn}[Class $\fclasscomplete{\aleph}{\mathnormal{\Omega}_z, \mathnormal{\Omega}_t}$]
	A function $(z,t) \mapsto$ \\ $f(z, t)$ is an element of the function class $\fclasscomplete{\aleph}{\mathnormal{\Omega}_z, \mathnormal{\Omega}_t}$
	if $f(\cdot, t) \in \mathcal{C}^\infty(\mathnormal{\Omega}_z)$ and $f(z, \cdot) \in \gevclasscomplete{\aleph}{\mathnormal{\Omega}_t}$.
\end{defn}

\section{State Feedback}
\label{sec:backstepping}
The following section states the main result of this contribution concerning the backstepping-based
state feedback.

\subsection{Error System}
	For the system \eqref{eq:transformed_system}, let the error coordinates be given by 
	\begin{subequations}
		\begin{align}
			\label{eq:error_temp}
			\tterrgz &= \ttransgz - \ttransrefgz \\
			\label{eq:error_pos}
			\dzi &= \zi - \zir,
		\end{align}
	\end{subequations}
	yielding the nonlinear error dynamics
		\begin{align}
			\partial_t \tterrgz 
			&= \hdg\partial_{\zt}^2 \tterrgz 
			+ \big( \dvi\partial_{\zt} + \vir\partial_{\zt} \big) \tterrgz \nonumber\\
			&\quad+ \dvi\partial_{\zt} \ttransrefgz
			\label{eq:error_pde}%
		\end{align}
	as shown in Appendix \ref{app:error_sys}. 
	Linearising \eqref{eq:error_pde} around the reference $(\ttnerrgz \equiv 0, \dzi=0)$ yields
	the linearised error dynamics
	\begin{align}
		\label{eq:approx_error_pde}%
		\partial_t \tterrgz &= 
		\hdg\partial_{\zt}^2 \tterrgz 
		+ \vir\partial_{\zt} \tterrgz \nonumber\\
		&\quad+ \dvi\partial_{\zt} \ttransrefgz
	\end{align}
	as well as
	\begin{equation}
		\dvi = 
		\scg \partial_{\zt}\tterrg(0, t) 
		+ \scag \partial_{\zt}\tterrag(0, t). 
		\label{eq:error_sc}
	\end{equation}
	Hence, substituting \eqref{eq:error_sc} in \eqref{eq:approx_error_pde} gives 
	the linear time-variant error system
	\begin{subequations}
		\begin{align}
			\label{eq:lin_error_pde}%
			\partial_t \tterrgz
			&= \hdg\partial_{\zt}^2\tterrgz 
			+ \vir\partial_{\zt} \tterrgz \nonumber\\
			&\quad+ b_{\gen}(\zt, t) \partial_{\zt}\tterrg(0,t)
			+ c_{\gen}(\zt, t) \partial_{\zt}\tterrag(0,t) \\
			\mkern-100mu%
			\label{eq:lin_error_bc_1}%
			\partial_{\zt} \tterrg(\zbtg,t) &= \inptg \\
			\label{eq:lin_error_bc_2}%
			\tterrg(0,t) &= 0 \\
			\dvi &= 
			\scg \partial_{\zt}\tterrg(0, t) 
			+ \scag \partial_{\zt}\tterrag(0, t) 
		\end{align}
		\label{eq:lin_error_sys}%
	\end{subequations}
	with 
	$b_{\gen}(\zt, t) =\scg\partial_{\zt} \ttransrefgz$,
	$c_{\gen}(\zt, t) = \scag\partial_{\zt} \ttransrefgz$, and the new input $\inptg$.
	Note that herein, the solid and liquid temperature errors are coupled over the whole domain by
	means of their fluxes through the phase boundary.

\subsection{Hopf-Cole Transformation}
	Since the coefficient $\vir$ of the convection term cannot be treated via the classic 
	backstepping transform,
	it is eliminated by a Hopf-Cole transformation \citep{Hopf1950}
	\begin{equation}
		\tterrgz = \hctransgz\therrgz
		\label{eq:hc_trafo}
	\end{equation}
	and choosing $\hctransgz = \exp\left(-\frac{\vir}{2\hdg}\zt\right)$, 
	which is a standard procedure in these cases.
	This yields the system
	\begin{subequations}
		\begin{align}
			\label{eq:hc_error_pde}%
			\mkern-20mu\partial_t \therrgz
			&= \hdg\partial_{\zt}^2\therrgz 
			+ r_{\gen}(\zt, t) \therrgz \nonumber \\
			&\;\;\; + \bar{b}_{\gen}(\zt, t) \partial_{\zt}\therrg(0,t)
			+ \bar{c}_{\gen}(\zt, t) \partial_{\zt}\therrag(0,t)
			\\
			\label{eq:hc_error_bc_1}%
			\partial_{\zt}\therrg(\zbtg, t) &=
			\frac{\vir}{2\hdg} \therrg(\zbtg, t)
			+ \hctransig(\zbtg,t)\inptg
			\\
			\label{eq:hc_error_bc_2}%
			\therrg(0,t) &= 0 \\
			\label{eq:hc_error_sc}%
			\dvi &=\scg \partial_{\zt}\therrg(0, t) 
			+ \scag \partial_{\zt}\therrag(0, t) 
		\end{align}%
		\label{eq:hc_error_sys}%
	\end{subequations}
	where
	\begin{subequations}
		\label{eq:hc_variables}
		\begin{align}
			\label{eq:r}
			r_0(\zt, t) = -\frac{1}{4\hdg} \left(
				2\aicomplete{r}{}{t}\zt + \vicomplete{r}{2}{t}
			\right) \\ 
			\label{eq:bbar}
			\bar b_{\gen}(\zt, t) = \hctransigz b_{\gen}(\zt,t)\\
			\label{eq:cbar}
			\bar c_{\gen}(\zt, t) = \hctransiagz c_{\gen}(\zt, t).
		\end{align}
	\end{subequations}
	Note, that the resulting system now exhibits reactive terms that are driven by the reference interface velocity and acceleration.
	Furthermore, the original Neumann boundary condition now appears as a Robin boundary condition.

\subsection{Backstepping Transformation}
	To enforce proper tracking of the reference, the errors in temperature and boundary position
	should converge to zero. 
	This demand is formulated in the target system
	\begin{subequations}
		\begin{align}
			\label{eq:tar_pde}
			\partial_t \ttargz
			&= \hdg\partial_{\zt}^2\ttargz 
			+ \mu_{\gen}(\zt, t) \ttargz  \\
			\label{eq:tar_bc_1}
			\partial_{\zt} \ttarg(\zbtg,t) &=  \bfd\nu_{\gen}\ttarg(\zbtg, t)\\
			\label{eq:tar_bc_2}
			\ttarg(0,t) &= 0 
		\end{align}
		\label{eq:tar_sys}%
	\end{subequations}
	with the reaction coefficient $\mu_{\gen}(\zt,t)$ and boundary gain $\nu_{\gen}$
	as design parameters.
	To map the system \eqref{eq:hc_error_sys} into \eqref{eq:tar_sys} the transformation
	\begin{equation}
		\label{eq:backstepping_trafo}
		\ttargz = \therrgz - \int
		_0^{\zt} 
        \ktgzzt \therrg(\zeta, t) \dop{\zeta}
	\end{equation}
	is used.
	Computing the requirements on the transformation kernel $\tilde k_{\gen}(\zt, \zeta,t)$
	(cf.~Appendix \ref{app:kernel_eqs}) yields the kernel equations 
	\begin{subequations}
		\begin{align}
			\label{eq:kernel_pde}
			\partial_t \ktgzzt
			&= \hdg \left(\partial_z^2 \ktgzzt 
			- \partial_{\zeta}^2 \ktgzzt\right)\nonumber
			\\
			&\quad +a_{\gen}(\zeta, t) \ktgzzt \\
			\label{eq:kernel_bc_zz}
			2\hdg
			\kernelcomplete{\tilde k_{\gen}}{\zt}{\zt}{t} &= 
				\int
				_0^{\zt} 
				a_{\gen}(\zeta, t)
				\dop{\zeta} -2\bar b_{\gen}(\zt, t)
			\\
			\label{eq:kernel_bc_z0}
			\hdg
			\kernelcomplete{\tilde k_{\gen}}{\zt}{0}{t} &=
			\int
			_0^{\zt} \bar b_{\gen}(\zeta, t) 
				\ktgzzt
				\dop{\zeta} 
				-\bar b_{\gen}(\zt, t)
			\\
			\label{eq:kernel_bc_extra}
			0 &=
			\int
			_0^{\zt} \bar c_{\gen}(\zeta, t) 
				\ktgzzt
				\dop{\zeta} 
				-\bar c_{\gen}(\zt, t)
		\end{align}%
		\label{eq:kernel_sys}%
	\end{subequations}
	where $a_{\gen}(\zt, t) = \mu_{\gen}(\zt, t) - r_{\gen}(\zt, t)$.
	Examining the system \eqref{eq:kernel_sys}, one observes that the problem for 
	$\ktgzzt$, given by 
	\eqref{eq:kernel_pde}--\eqref{eq:kernel_bc_z0} is well-posed (cf.~\ref{ssec:existance}).
	Herein, the integral boundary condition \eqref{eq:kernel_bc_z0} arises since the term 
	$\bar b_{\gen}(\zt,t)\partial_{\zt}\therrg(0,t)$ is to be eliminated from \eqref{eq:hc_error_pde}.
	However, the demand for completely decoupled target systems and, thus, the elimination of
	$\bar c_{\gen}(\zt,t)\partial_{\zt}\therrag(0,t)$ that results in \eqref{eq:kernel_bc_extra}
	renders the problem overdetermined.
	Therefore, to recover a well-posed formulation the convective coupling at $\zt=0$ is reintroduced
	with the modified target system dynamics
	\begin{multline}
		\label{eq:mtar_pde}
		\partial_t \ttarmgz
		= \hdg\partial_{\zt}^2\ttarmgz \\
		+ \mu_{\gen}(\zt, t) \ttarmgz
		+ d_{\gen}(\zt, t) \partial_{\zt}\ttarmag(0,t)
	\end{multline}
	where
	$d_{\gen}(\zt,t) =  \int^{\zt} \bar c_{\gen}(\zeta, t) 
		\ktgzzt \dop{\zeta} 
		-\bar c_{\gen}(\zt, t)$
	replaces \eqref{eq:kernel_bc_extra}.
	Obviously, by choosing $\mu_{\gen}(\zt,t) \le 0 \; \forall (\zt, t)$ and $\nu \le 0$ this 
	approach yields exponentially stable error dynamics for the one-phase case
	where the gradient in the adjacent phase vanishes from \eqref{eq:error_sc}.
	For the two-phase case, stability of the resulting error dynamics has to be shown due to the bilateral
	coupling via $d_{\gen}(\zt,t)$.
	This will be addressed in a forthcoming publication due to lack of space.
	However, simulation studies yield promising results as Section \ref{sec:results} shows.
	Certainly, for both variants the target systems properties can only be conveyed if the inverse
	transformation of \eqref{eq:backstepping_trafo} exists.
	This can be assumed since it is of Volterra-type and therefore always invertible 
	(cf.~\citealt{Heuser1992}) or shown by a simple fixed-point argument\footnote{Note, that to rigorously proof the invertability, the transformation has to be defined as a linear map on an appropriately chosen Banach space in order to make the underlying fixed point theory applicable. For lack of space, details are omitted within this contribution.}.

	Finally, by examining \eqref{eq:tar_bc_1}, eliminating the target terms via 
	\eqref{eq:backstepping_trafo} and substituting \eqref{eq:hc_error_bc_1},
	the control input for the original system 
	with the kernel in original coordinates 
	$\kernelcomplete{k_{\gen}}{z}{\zeta}{t} = \kernelcomplete{\tilde k_{\gen}}{\tilde z}{\zeta}{t}$
	is given by 
	\begin{align}
	\label{eq:orig_input}
		\inpg &= \frac{\hcg}{\bfd}\Bigg[
			\left( \kernelcomplete{k_{\gen}}{\zbg}{\zbg}{t} 
				+ \bfd\nu_{\gen}
				- \frac{\vir}{2\hdg} 
			\right)\times \nonumber\\
			&\qquad\qquad\left(
				\vphantom{T^2}
			\torigg(\zbg,t) - \trefg(\zbg,t)\right)
			+ \partial_{z}\trefg(\zbg,t)
			\nonumber\\ 
			&\qquad
			+ \int
			_{\zi}^{\zbg} \left(
				\vphantom{T^2}
				\partial_{z}\kernelcomplete{k_{\gen}}{\zbg}{\zeta}{t} 
				+ \bfd\nu\kernelcomplete{k_{\gen}}{\zbg}{\zeta}{t} 
			\right)
			\times\nonumber\\ &\mkern-18mu
			\left(
				\vphantom{T^2}
			\torig(\zeta,t) - \trefg(\zeta,t)\right)
			\exp\left(\frac{\vir}{2\hdg}\left(\zeta-\zbg\right)\right)
			\dop{\zeta}
		\Bigg].
	\end{align}

\section{Well-posedness and numerical solution of the kernel equations}
\label{sec:kernel}
	In this section, a numerical solution scheme for the kernel equations 
    \eqref{eq:kernel_pde}--\eqref{eq:kernel_bc_z0} is discussed\footnote{
    As only variables of one phase occur in the kernel equations, the generic placeholder $\gen$ is dropped
    in favour of a more compact notation from now on.
    }.
	To this end, the integral form of the system is derived in a first step.
	As stated in \citep{Jadachowski2012}, the method of successive approximations as introduced in \citep{Colton1977} and extended in \citep{Meurer2009} does not
	show good convergence for time-varying kernels and is therefore only employed to investigate the
	existence of a solution.
	Thus, the presented solution will be based on a spatial discretisation of the kernel.
	However, in contrast to \citep{Jadachowski2012} the Midpoint rule will be used which eventually leads
	to an iterative solution scheme that maintains the structural properties of the problem.

	\subsection{Integral Form}
	Introducing the normal form coordinates
	$\eta = \zt + \zeta$ and $\sigma = \zt - \zeta$
	yields the dynamics
	\begin{subequations}
		\begin{align}
			\label{eq:kernel_pde_4}
			\partial_t \kbest 
			\MoveEqLeft[1]
            = 4\hd\partial_{\eta\sigma}\kbest
			+ a\left(\tfrac{\eta-\sigma}{2}, t\right) \kbest \\
			\label{eq:kernel_bc_e0_4}
			4\hd
			\kbeot =& 
				\int
				_0^\eta 
                a\left(\tfrac{r}{2}, t\right)\dop{r}
				- 4\bar b\left(0, t\right)
			\\
			\label{eq:kernel_bc_ee_4}
			\hd
			\kbeet =& 
				\int
				_0^\eta \bar b\left(\tfrac{\eta+s}{2}, t\right)
				\kernelcomplete{\kb}{\eta}{s}{t}\dop{s}
				-\bar b\left(\eta,t\right)
		\end{align}
	\end{subequations}
	of the transformed kernel $\kbest = \ktzzt$.
	Hence, formal integration of \eqref{eq:kernel_pde_4} wrt.~$\sigma$ and $\eta$ as well as substitution
	of \eqref{eq:kernel_bc_e0_4} (derived wrt.~$\eta$) and \eqref{eq:kernel_bc_ee_4} yields
	\begin{align}
		\label{eq:int_form_2}
		\kbest &=
		\frac{1}{4\hd}\left\{
			\int
			_{\sigma}^{\eta}
			\bigg[
			a\left(\tfrac{r}{2}, t\right)
			+\int
			_0^\sigma 
			\bigg(
			\partial_t \kernelcomplete{\kb}{r}{s}{t}
		\right. \nonumber\\
		&\hphantom{=\frac{4}\quad{\hd}\Bigg(}
			- a\left(\tfrac{r-s}{2}, t\right)
			\kernelcomplete{\kb}{r}{s}{t}{} 
			\bigg) \dop{s}
			\bigg] \dop{r}
			\left. \vphantom{\int_{\sigma}^{\eta}} \right\}
		\nonumber\\
		&\mkern-60mu + \frac{1}{\hd}\left(
			\int
			_0^\sigma \bar b\left(\tfrac{\sigma+s}{2}, t\right)
			\kernelcomplete{\kb}{\sigma}{s}{t}\dop{s}
			- \bar b\left(\sigma,t\right)
		\right).
	\end{align}%

	\subsection{Existence}
	\label{ssec:existance}
	As already presented in \citep{Colton1977} a solution for \eqref{eq:int_form_2} by means of the method of successive approximations can be established by considering the series
	\begin{equation}
		\kbest = \sum\limits_{n=0}^{\infty} \Kbest{n}
		\label{eq:approx_series}
	\end{equation}
	where the $\Kbest{n}$ are given by $\Kbest{0} = 0$,
	\begin{subequations}
		\begin{align}
			\Kbest{1} & = 
			\int
			_{\sigma}^{\eta} a\left(\tfrac{r}{2}, t\right) \dop{r}
			- \frac{1}{\hd} \bar b\left(\sigma,t\right)
			\\
			\label{eq:kn}
			\Kbest{n} &=
			\frac{1}{\hd} \int
			_0^\sigma 
				\bar b\left(\tfrac{\sigma+s}{2}, t\right)
				\Kernelcomplete{\Kb}{\sigma}{s}{t}{n-1}\dop{s}\nonumber\\
			&\quad
			+\frac{1}{4\hd}\int
			_{\sigma}^{\eta}\int
			_0^\sigma 
			\left( \vphantom{\bigg)}
				\partial_t \Kernelcomplete{\kb}{r}{s}{t}{n-1}
			\nonumber\right.\\ &
			\qquad
				- a\left(\tfrac{r-s}{2}, t\right)
				\Kernelcomplete{\Kb}{r}{s}{t}{n-1}
			\left.\vphantom{\bigg)}\right)\dop{s} \dop{r}
		\end{align}
	\end{subequations}
	Furthermore, \cite{Meurer2009} show that the convergence conditions on the series 
	\eqref{eq:approx_series} 
	depend on an upper bound for $\partial_t^l\Kbest{n}$ due to the repetitive
	differentiation in \eqref{eq:kn}.
	Although the detailed appearance may differ from the systems discussed in \citep{Meurer2009} or
	\citep{Izadi2015}, the terms that are to be examined are $a(\zt,t)$ and
	$\bar b(\zt,t)$.

	Consider the term $a(\zt, t)$, which is composed by the design parameter $\mu(\zt, t)$
	and the reaction coefficient $r(\zt, t)$ of the error system.
	The latter coefficient, however, is mostly characterised by the reference interface velocity and
	acceleration which --as design parameters-- are characterised by 
	$\zir \in \gevclasscomplete{\aleph\le2}{\mathbb{R}^+}$, (cf. \ref{eq:r}).
	Thus, by demanding 
	$\mu(\zt,t) \in \fclasscomplete{\aleph\le2}{\Xi}$ 
	with $\mathnormal{\Xi} = (0, \zbtg) \times \mathbb{R}^+$ it follows that 
	$a(\zt,t) \in \fclasscomplete{\aleph\le 2}{\mathnormal{\Xi}}$.
	Next, $\bar b(\zt,t)$ from \eqref{eq:bbar} is given by the product of the reference
	temperature gradient
	$\partial_{\zt}\ttransrefgz \in \fclasscomplete{\aleph\le2}{\mathnormal{\Xi}}$
	and $\hctransigz$.
	The inverse of the transformation \eqref{eq:hc_trafo} in turn is a simple composition of 
	$\vir$ and the exponential function which does not affect the Gevrey order \citep{Gevrey1918}.
	Hence, $\bar b(\zt,t)$ is from $\fclasscomplete{\aleph\le2}{\mathnormal{\Xi}}$.

	Summarising, the method employed in \citep{Meurer2009} can be adapted to this case to show existence for
	and uniqueness of the solution.

	\subsection{Approximation Scheme}
	For the spatial discretisation, values of the kernel are computed at selected grid points, given 
	by $\bar k_{i,j}(t) = \kernelcomplete{\kb}{i\discre}{j\discrs}{t}$ with the grid widths $\discre$ and
	$\discrs$ in $\eta$- and $\sigma$-direction, respectively.
	Due to the varying extent of the kernel domain, the computation grid is taken as the maximum extent
	for each respective phase, given by $\Delta\Gamma = \zbl- \zbs$.
	This yields the number of nodes $N_{\eta}=2\Delta\Gamma / \discre$ and  $N_{\sigma}=\Delta\Gamma / \discrs$ in each direction.
	\begin{figure}
		\centering
		\resizebox{\linewidth}{!}{\begin{tikzpicture}[auto, >=latex] 
	\pgfmathsetmacro{\dz}{4}

	\draw[ultra thick, color=HKS92K10, fill=HKS92K10, opacity=1] (0, 0) 
		-- (1.25*\dz, 0)
		-- (1.25*\dz, .75*\dz)
		-- (.75*\dz, .75*\dz)
		-- (0, 0);
	\draw[ultra thick, color=HKS92K30,fill=HKS92K30, opacity=1] (0, 0) 
		-- (\dz, 0)
		-- (\dz, .5*\dz)
		-- (.5*\dz, 0.5*\dz)
		-- (0, 0);
	\draw[ultra thick, color=HKS92K50, fill=HKS92K50, opacity=1] (0, 0) 
		-- (.75*\dz, 0)
		-- (.75*\dz, .25*\dz)
		-- (.25*\dz, .25*\dz)
		-- (0, 0);
	\node[] at (1.125*\dz, .6125*\dz) {$d=0$};
	\node[] at (0.8625*\dz, .375*\dz) {$d=1$};
	\node[] at (0.6125*\dz, .125*\dz) {$d=2$};

        \draw[->] (0,0) -- (2.20*\dz,0) node[below] {$\eta$};
        \draw[->] (0,0) -- (0,1.15*\dz) node[left] {$\sigma$};

        \draw[-, thick] (0, 0) -- (\dz,\dz) -- (2*\dz, 0) -- (0, 0);

        \draw[-, ultra thick, color=HKS44K100] (0, 0) -- node [above left] {$\sigma=\eta$} (\dz,\dz);
        \draw[-, ultra thick, color=HKS07K100] (0, 0) -- node [below] {$\sigma=0$}(2*\dz, 0);

	\draw[<->] (1.5*\dz, .2*\dz) -- node [below] {$\discre$} (1.75*\dz,.2*\dz);
        \draw[<->] (1.3*\dz, 0*\dz) -- node [right] {$\discrs$} (1.3*\dz,.25*\dz);

	\foreach \i in {0.0, .25, .5, ..., 2}{
		\tikzmath{\a = \i * \dz;}
		\tikzmath{\b = \i * \dz;}
		\newdimen\id
		\id = \i mm
		\ifdim \id > 1mm \tikzmath{let \b = 2*\dz - \i*\dz;} \fi
		\draw[dashed, color=HKS92K80] (\a,0) -- (\a, \b);
	}
	\foreach \i in {.25, .5, ..., 1}{
		\tikzmath{
			let \a1 = \i * \dz;
			let \a2 = 2*\dz - \i*\dz;
			let \b = \i * \dz;
		}
		\draw[dashed, color=HKS92K100] (\a1, \b) -- (\a2, \b);
	}

	\foreach \j in {0, ..., 4}{
		\tikzmath{\istart = \j;}
		\tikzmath{\iend = 8 -\j;}
		\foreach \i in {\istart, ..., \iend}{
			\tikzmath{\a = \i * \dz/4.0;}
			\tikzmath{\b = \j * \dz/4.0;}
			\fill[fill=black] (\a, \b) circle [radius=1.5pt] node [below right] {};
		}
	}


	\fill (1.25*\dz, .75*\dz) circle [radius=2pt] node [right] {$\bar k_{i, j}(t)$};

	\fill (0, 0) circle [radius=2pt] node [below ] {$\bar k_{0,0}(t)$};
	\fill (\dz, \dz) circle [radius=2pt] node [right] {$\bar k_{N_{\sigma}, N_{\sigma}}(t)$};
	\fill (2*\dz, 0) circle [radius=2pt] node [below] {$\bar k_{N_{\eta}, 0}(t)$};

\end{tikzpicture}}
		\caption{Discretised kernel domain with selected kernel elements and required  derivative orders for the computation of $k_{i,j}(t)$ (shaded).}
		\label{fig:spat_disc}
	\end{figure}
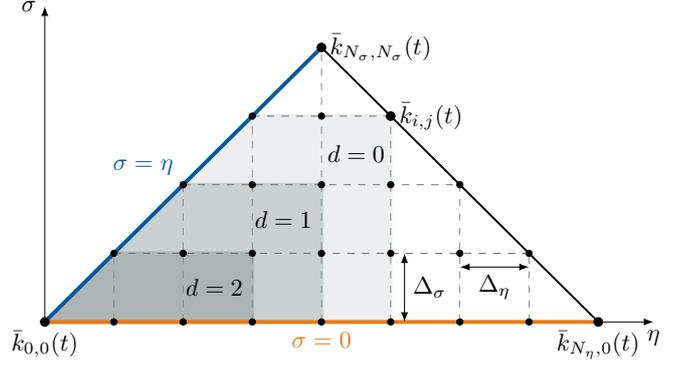

	Next, the uniform step width $\discr = \discre = \discrs$, motivated
	by the similar dynamics in both directions, is chosen.
	Hence, approximating the integrals in \eqref{eq:int_form_2} by lower sums
	yields the explicit computation scheme
	\begin{align}
		\label{eq:su_form}
		\kij{i}{j} &=
		\frac{\discr^2}{4\hd}
		\sum\limits_{n=j}^{i-1}
		\sum\limits_{m=0}^{j-1} \bigg(
			\mkern-2mu
			\partial_t \kij{n}{m}
			- a\left(\tfrac{n-m}{2}\discr, t\right)
			\mkern-2mu
			\kij{n}{m}
			\mkern-2mu
		\bigg) 
		\nonumber\\ &\quad
		+\frac{\discr}{4\hd} \sum\limits_{n=j}^{i-1}
			a\left(\tfrac{n}{2}\discr, t\right)
		- \frac{1}{\hd}\bar b\left(j\discr,t\right)
		\nonumber\\ &\quad
		+\frac{\discr}{\hd} \sum\limits_{m=0}^{j-1}
			\bar b\left(\tfrac{j+m}{2}\discr, t\right)
			\kij{j}{m}  
	\end{align}
    for the interior where $0<j<N_{\sigma}$ and $j<i<N_{\eta}-j$
	while the boundary expressions are given by
	\begin{equation}
		\label{eq:ki0}
		\kij{i}{0} 
        = \frac{\discr}{4\hd}\left(
			\sum\limits_{m=0}^{i-1}
			a\left(\tfrac{m}{2}\discr, t\right)
			- 4\bar b\left(0,t\right)\right)
	\end{equation}
	\begin{equation}
		\label{eq:kjj}
        \mkern-8mu
		\kij{j}{j} 
        = \frac{1}{\hd}
			\mkern-1mu
        \left(
			\mkern-1mu
			\discr
			\mkern-1mu
            \sum\limits_{n=0}^{j-1}
			\bar b\left(\tfrac{j+n}{2}\discr, t\right)
			\kij{j}{n}  
            - \bar b\left(j\discr,t\right)
        \right)
	\end{equation}
    for $0\le i \le N_{\eta}$ and $0\le j \le N_{\sigma}$, respectively.
	For the computation of an arbitrary kernel element $\kij{i}{j}$ via 
	\eqref{eq:su_form} the temporal derivatives of the neighbouring elements 
	still pose a problem.
	However, all temporal derivatives of $\kij{i}{j}$ can be recursively substituted
	until only derivatives of $a(\zt,t)$ and $\bar b(\zt,t)$ remain. 
	This gives the map 
	\begin{equation}
		\kij{i}{j} = \mathnormal{\Theta}_{i,j}\left(a^{(0)}, \dotsc, a^{(j)},
			\bar b^{(0)}, \dotsc, \bar b^{(j)}
		\right)
		\label{eq:ki_map}
	\end{equation}
	which is further discussed in Appendix \ref{app:derivatives}.

	Hence, given the structures of $a(\zt,t)$ and $\bar b(\zt,t)$, the complete kernel 
	can be expressed as a function of the reference trajectory for the flat output $\y$ of 
	\eqref{eq:transformed_system}, its derivatives, and the design parameter $\mu(\zt,t)$.
	Note that in \citep{Jadachowski2012} this property is lost since an initial value problem
    has to be solved for the inner kernel elements,
	while the presented approach converges to the solution via successive approximations
	\eqref{eq:approx_series} for sufficiently small step sizes.
	However, the usage of the trapezoidal rule as in \citep{Jadachowski2012}
	drastically reduces
	the approximation error for similar grid sizes due to the implicit nature of the resulting
	approximation scheme.

\section{Simulation Setup and Results}
\label{sec:results}

This section briefly presents simulation results for the two-phase case.
The simulated process goal is the growth of a \gls{gaas} single crystal in a \SI{400}{\milli\metre}
furnace.
As simulation model, a lumped, FE-based approximation with $41$ nodes per phase has been used.
For further details on the physical parameters of the simulation setup, 
please refer to \cite[Sec.~7]{Ecklebe2019}.

In this case, the growth process shall start at an initial length of \SI{200}{\milli\metre}
and end at about \SI{300}{\milli\metre}.
To do so, $y_{\mathrm{r}, 2}(t)$ 
is chosen as a smooth transition (cf.~\eqref{eq:flat_traj})
between these values with the targeted solid phase gradient at the interface held constant at 
$y_{\mathrm{r}, 1}(t) = \SI{17}{\kelvin\per\meter}$ during the transition ($t < \SI{25}{\hour}$). 
Thus, the recursion \eqref{eq:coeff_iteration} from Section \ref{sec:feedforward} yields the
reference temperature trajectory $t \mapsto \tref(\cdot,t)$ which in turn enables the
computation of the kernels for the solid and liquid phases by means of \eqref{eq:ki_map}.
Herein, \SI{81} discrete points and the parameters $\mu(\zt,t)=\SI{-1e-2}{\per\second}$ 
as well as $\nu=\SI{0}{\per\meter}$ have been used in the target systems of both phases
(cf.~\eqref{eq:tar_sys}).
Finally, the control input is computed via \eqref{eq:orig_input}.

\begin{figure}
	\centering
	\resizebox {\linewidth} {!} {
		\input{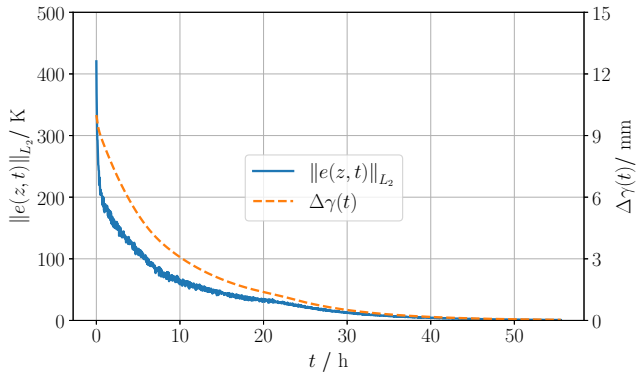}
	}
	\caption{Trajectories of the temperature errors $L_{2}$ norm (solid)
	and the interface deviation $\dzi$ (dashed).}
	\label{fig:errors}
\end{figure}

For demonstration purposes, an initial error for the phase boundary of 
$\Delta\gamma(0) = \SI{10}{\milli\meter}$ as well $\Delta \dot \gamma (0) = \SI{-3}{\milli\meter\per\hour}$ 
for the growth velocity have been introduced.
As shown in Figure \ref{fig:errors}, the controller ensures convergence for the error system
$\terrz$ as well as the interface deviation $\dzi$.

\section{Conclusion and outlook}

In this contribution, a backstepping-based tracking control has been presented for a two-phase Stefan problem,
occurring in the \gls{vgf} process. 
By utilising both system inputs, this approach enables the tracking of arbitrary\footnote{
The only hard restriction on the chosen references is given by the Gevrey order condition which is rather easy to satisfy.}reference
trajectories.
Although only a linear approximation of the error system is used, the results are promising.
However, the stability proof of the modified target system for the two-phase case remains open.
Moreover, for output feedback in practical applications, the proposed feedback has to be
complemented by an appropriate observer.
While two early-lumping observers are presented in \citep{Ecklebe2019} by the authors for the \gls{tpsp}
and a late-lumping design is given in \citep{Koga2019a} for the \gls{opsp},
the treatment of the two-phase case with the backstepping method still remains open.
These issues are currently under investigation and will be addressed in a more detailed article.


\bibliography{bibliography}
\appendix
\section{Error System dynamics}
\label{app:error_sys}
Taking the temporal derivative of \eqref{eq:error_temp} and substituting \eqref{eq:transformed_pde} for
the original and reference temperature yields
\begin{align}
	\partial_t \tterrgz 
	&= \partial_t \ttransgz - \partial_t\ttransrefgz \nonumber\\
	&= \hdg\partial_{\zt}^2 \ttransgz + \vi\partial_{\zt} \ttransgz\nonumber \\
	&\quad- \hdg\partial_{\zt}^2\ttransrefgz - \vir\partial_{\zt} \ttransrefgz \nonumber
	\intertext{which by substituting \eqref{eq:error_temp}, derived wrt.~$\zt$ twice, becomes:}
	&
	\mkern-18mu
	= \hdg\partial_{\zt}^2 \tterrgz 
		+ \vi\partial_{\zt} \ttransgz 
		- \vir\partial_{\zt} \ttransrefgz \nonumber.
	\intertext{Furthermore, replacing $\vi$ and $\ttransgz$ via \eqref{eq:error_temp} and \eqref{eq:error_pos}, respectively, gives}
	\partial_t \tterrgz 
	&= \hdg\partial_{\zt}^2 \tterrgz 
		- \vir\partial_{\zt} \ttransrefgz \nonumber\\
		&\quad+ \Big(\dvi + \vir\Big)\left(\partial_{\zt} \tterrgz + \partial_{\zt} \ttransrefgz\right)\nonumber
\intertext{and thus}
	\partial_t \tterrgz 
	&= \hdg\partial_{\zt}^2 \tterrgz 
		+ \big( \dvi\partial_{\zt} + \vir\partial_{\zt} \big) \tterrgz\nonumber \\
		&\quad+ \dvi\partial_{\zt} \ttransrefgz.
\end{align}
Finally, examining the boundaries yields
\begin{align}
	\label{eq:error_bc_1}
	\tterrg(0,t) &= \ttransg(0,t) - \ttransrefg(0, t) = \tm - \tm = 0\\
	\label{eq:error_bc_2}
	\partial_{\zt}\tterrg\left(\zbtg,t\right) &=
	\partial_{\zt}\ttransg\left(\zbtg,t\right) - \partial_{\zt}\ttransrefg\left(\zbtg,t\right)\nonumber\\
	&= \frac{\bfd}{\hcg}\inpg -\partial_{\zt}\ttransrefg(\zbt,t) =\inptg.
\end{align}

\section{Kernel Equations}
\label{app:kernel_eqs}
	To derive the mandatory conditions on the transformation kernel from \eqref{eq:backstepping_trafo},
	the derivatives of \eqref{eq:backstepping_trafo} are taken wrt.~$t$ and $\zt$.
	Next, they are substituted in \eqref{eq:tar_pde} to express the target system dynamics in terms of the 
	original error $\tterrgz$.
	Furthermore, \eqref{eq:error_pde} is substituted and integration by parts is performed to shift the
	spatial operators onto $\tilde k(\zt, \zeta, t)$.
	Finally, since the resulting equation has to hold for all $\tterrgz$, one arrives at the kernel equations
	\begin{subequations}
		\begin{align}
			\label{eq:kernel_pde_0}
			\partial_t \kernelcomplete{k_{\gen}}{\zt}{\zeta}{t}
			&= \hdg \left(\partial_z^2 \kernelcomplete{k_{\gen}}{\zt}{\zeta}{t} 
			- \partial_{\zeta}^2 \kernelcomplete{k_{\gen}}{\zt}{\zeta}{t}\right)\nonumber
			\MoveEqLeft[10]
			\\
			&\quad +a_{\gen}(\zeta, t) 
				\kernelcomplete{k_{\gen}}{\zt}{\zeta}{t} \\
			\label{eq:kernel_bc_zz_0}
			a_{\gen(\zt,t)}
			&=
			2\hdg
			\tfrac{\textup{d}}{\textup{d}\zt} \kernelcomplete{k_{\gen}}{\zt}{\zt}{t} 
			\\
			\label{eq:kernel_bc_z0_0}
			\hdg
			\kernelcomplete{k_{\gen}}{\zt}{0}{t} &=
			\int_0^{\zt} \bar b_{\gen}(\zeta, t) 
				\kernelcomplete{k_{\gen}}{\zt}{\zeta}{t}
				\dop{\zeta} 
				-\bar b_{\gen}(\zt, t)
			\\
			\label{eq:kernel_bc_extra_0}
			0 &= \int_0^{\zt} \bar c_{\gen}(\zeta, t) 
				\kernelcomplete{k_{\gen}}{\zt}{\zeta}{t}
				\dop{\zeta} 
				-\bar c_{\gen}(\zt, t).
		\end{align}%
	\end{subequations}%
	Thus, integration of \eqref{eq:kernel_bc_zz_0} yields
	\begin{equation*}
		\kernelcomplete{k_{\gen}}{\zt}{\zt}{t} = 
		\frac{1}{2\hdg} \int_0^{\zt} a_{\gen}(\zeta, t) \dop{\zeta}
		+ \kernelcomplete{k_{\gen}}{0}{0}{t}
	\end{equation*}
	where $\kernelcomplete{k_{\gen}}{0}{0}{t} = -b_{\gen}(\zt, t)/\hdg$ from
	\eqref{eq:kernel_bc_z0_0} at $\zt=0$.

\section{Derivative orders}
\label{app:derivatives}
	To analyse the required derivative orders for the computation of an arbitrary kernel element
	$\kij{i}{j}$ the $l$-th derivative of \eqref{eq:su_form}
	\begin{align}
		\partial_t^l\kij{i}{j} &=
		\frac{\discr^2}{4\hd}
		\sum\limits_{n=j}^{i-1}
		\sum\limits_{m=0}^{j-1} \Bigg(
		\partial_t^{(l+1)} \kij{n}{m} \nonumber\\
		&\mkern-36mu
		-\sum\limits_{d=0}^{l}\binom{l}{d} 
		\partial_t^{(l)}a\left(\tfrac{n-m}{2}\discr, t\right)
		\partial_t^{(l-d)}\kij{n}{m} \Bigg)\nonumber\\
		&\mkern-36mu +\frac{\discr}{\hd}
		\sum\limits_{m=0}^{j-1}
		\sum\limits_{d=0}^{l}\binom{l}{d} 
		\partial_t^{(l)}\bar b\left(\tfrac{j+m}{2}\discr, t\right)
		\partial_t^{(l-d)}\kij{j}{m}
		\nonumber\\ &\mkern-36mu
		-\frac{1}{\hd}\partial_t^l\bar b\left(j\discr,t\right)
		+\frac{\discr}{4\hd} 
		\sum\limits_{n=j}^{i-1}
		\partial_t^l a\left(\tfrac{n}{2}\discr, t\right)
	\end{align}
	is examined.
	As it can be seen, for every derivative in \eqref{eq:su_form} that is to be eliminated, 
	derivatives of the kernel elements below and left of one order higher, as well as derivatives
	of the provided functions $a(\zt,t)$ and $\bar b(\zt,t)$ with the same order are introduced
    (cf.~shaded areas in Figure \ref{fig:spat_disc}).
	Therefore, the computation of $\kij{i}{j}$ demands a derivative order of $\min \{i,j\} = j$ for
	these functions.
	Hence, for the computation of the complete kernel, $d_{\mathrm{max}} = N_{\sigma}-1$ derivatives
	are required as the kernel element $\kij{N_{\sigma}}{N_{\sigma}}$ is already given by \eqref{eq:kjj}.

\end{document}